\documentclass[12pt]{article}
\usepackage[T1,T2A]{fontenc} %
\usepackage[cp1251]{inputenc} %
\usepackage[russian]{babel} %
\usepackage{amsfonts} %
\usepackage{amsmath, bm} %
\usepackage{amssymb} %
\usepackage{MnSymbol}
\usepackage{color}
\usepackage[active]{srcltx}
\usepackage{xr}
\usepackage[final]{pdfpages}
\usepackage{graphicx}
\usepackage{wrapfig}
\usepackage{tikz}
\usetikzlibrary{arrows,decorations.pathmorphing,
backgrounds,positioning,fit,petri}
\usepgflibrary{arrows}
\usetikzlibrary{arrows}
\definecolor{whitesmoke}{rgb}{0.9,0.9,0.9}
\definecolor{whitesmoke2}{rgb}{0.8,0.8,0.8}
\usetikzlibrary{backgrounds}
\begin{document}


\smallskip	
\begin{center} {\large \bf On Uniqueness of a Solution to the Boundary Initial Value Problem} \end{center}
\begin{center} VLADIMIR V. BASOV \end{center}

\noindent
{\footnotesize
St.Petersburg State University, Universitetskaya nab., 7-9, St.Petersburg, 199034, Russian Federation}

\bigskip
{\small 
Abstract. \ A first-order ordinary differential equation, solved with respect to derivative, is considered. It's right-hand side is defined and continuous on the set, consisting of a connected open subset of a two-dimensional Euclidean space and a part of its boundary.
In the papers, dated by 2020, problems related to the existence or absence of the solution to the BIVP --- Initial Value Problem set at the boundary point --- were researched using different approaches. Such formulation of the Initial Value Problem is different from the formulation, established in the classic theory, where it is set at an interior point. This paper is devoted to solving problems related to uniqueness  or non-uniqueness of the BIVP solutions. New definitions, related  to uniqueness, absent in the formulation of the IIVP --- Initial Value Problem set at the interior point of the equation's domain --- are introduced. The theorems about the formal, local and global uniqueness of the BIVP solutions are proven. The differences between BIVP and IIVP are shown. For example, non-equivalence of the definitions of the formal, local and global uniqueness' for BIVP and IIVP is demonstrated, This non-equivalence leads to the appearance of hidden non-uniqueness points along with uniqueness and non-uniqueness points. Suggested theory is supposed to fill in the blanks in an existing literature, related to the problems of existence and uniqueness of the BIVP solutions.

}

\bigskip 
\begin{center}{\large 1. Introduction.} \end{center} 
\noindent\!1.1. {\bf The statement of the problem.} \
Consider the first-order ordinary differential equa- tion, solved with respect to derivative
$$y'=f(x,y),\eqno (1)$$
where $f(x,y)$ is the real function, continuous on set \,$\widetilde G=G\cup \widehat G,$\,
where $G$ is a connected open subset of the Euclidean space $\mathbb{R}^2,$
and the set $\widehat G$ is a subset of the boundary $\partial G$ of the open connected subset $G.$
Additionally, the <<good>> boundary $\widehat G$ contains all points of $\partial G,$ where $f(x,y)$ can be defined as continuous function. 

The set $\widetilde G$ is chosen as a domain of the equation (1), because in specific differential equations written with the help of elementary functions, the domain of the function $f$ does not have to be an open set and may contain the boundary points.
In such case, the solutions may both begin at the boundary and attain it.
At the same time, classical theorems 
guarantees existence of a solution that begins only at interior points.

The conditions for existence of the solution of the initial value problem set at the point $(x_0,y_0)\in \widehat G,$ and its non-existence are provided in the paper [1].

This can be achieved by extracting the cases, in which it is possible to construct the analogues of the Peano triangle and segment and apply the Euler polygonal method. The advantage of this method is its constructivity and simplicity.
\nopagebreak

However, the specific constraints of the Euler polygonal method do not allow to obtain theoretically more universal results.
Therefore, in the paper [2] the equation (1) was redefined in a such way, that the initial value problem, set at the boundary point $(x_0,y_0),$ became the interior initial value problem, which is a subject of the standard Peano theorem's application.

The so-called comparison theorems and differential inequalities [3,4] were used to solve the problem of the solution of the modified initial value problem being the solution of the initial boundary problem. This allowed to weaken the requirements on the boundary curves in the neighbourhood of the point $(x_0,y_0)$ and consider the cases, which could not be studied using the Euler polygonal method.

It should be noted that the research of the solutions of the boundary initial value problem is not considered in any academic or scientific books, despite its classical formulation. Thus, the papers [1] and [2] fill in the existing blanks.  

\smallskip
{\bf Remark 1.}\, 
{\sl The similar situation takes place in the paper [5],  in which the new results on the classical problem on the properties of the boundaries of the solution's maximal interval of existence, where the solutions are studied as the functions of initial values, are provided.   }

\smallskip
Along with the problems on existence or non-existence of the solutions of the boundary initial value problem, the problems about its uniqueness arise. Solutions to these problems can not be obtained applying the existing theorems, proven for the interior initial value problem. 
This paper is devoted to the study of these problems.

The first attempt to systematically present the new definitions and some results related to uniqueness at the boundary points, was provided in the chapter\;1 of the academic book [6]. This paper provides the partially revised and substantially updated theory, developed to solve the problems on uniqueness or non-uniqueness of the solutions of the boundary initial value problem. 

\medskip
\noindent\!1.2. {\bf On existence of the solution of the BIVP.} 
Let us formulate the necessary definitions and results from the paper [1]. 

\smallskip
{\bf Df.}\, 
Denote function $y=\varphi(x),$ defined on interval $\langle a,b\,\rangle ,$
as a {\it solution of the differential equation} $(1),$ if the following conditions hold: \
$1)$ the function $\varphi(x)$ is differentiable at any point $x\in\,\langle a,b\,\rangle ;$
$2)$ a point $(x,\varphi(x))\in \widetilde G$ for all $x\in\,\langle a,b\,\rangle ;$
$3)$ $\varphi'(x)=f(x,\varphi(x))$ for any $x\in\,\langle a,b\,\rangle $ \ 
(\,here symbol $\langle\,$ --- is $($ or $[\;,$ $\rangle\,$~--- analogously; $[a,b]$ is a segment, $(a,b)$ is an interval and $a<b).$

\smallskip
Now, the notion of a solution of the equation (1) can be specified depending on the location of its graph in the set $\widetilde G.$

\smallskip
{\bf Df.}\, 
Denote solution $y=\varphi(x)$ of the equation $(1),$ defined on interval $\langle a,b\,\rangle,$ as:
$a)$ {\it an interior solution,}\, if a point $(x,\varphi(x))\in G$ for any $x\in \langle a,b\,\rangle;$
$b)$~{\it a boundary solution,}\, if $(x,\varphi(x))\in \widehat G$ for any $x\in \langle a,b\,\rangle;$
$c)$~{\it a mixed solution,}\, if there are such $x_1,x_2\in\langle a,b\,\rangle,$ that $(x_1,\varphi(x_1))\in G$
and $(x_2,\varphi(x_2))\in \widehat G.$
	
\smallskip
Obviously, if the equation (1) is considered only in a connected open subset $G,$ then the definitions of a solution 
and an interior solution are the same.

\smallskip
{\bf Remark 2.}\, 
{\sl Usually functions $f(x,y),$ that are dealt with when solving a specific equations (1), are a composition of elementary functions that are continuous on their domain, which is often no more than countable union of the connected sets $\widetilde G.$ This is one of the reasons to consider the equation (1) on the set $\widetilde G$ instead of a connected open set.

For example, the equation $y'=\sqrt{y},$ which domain is the set $\widetilde G=\{(x,y)\!:\,x\in \mathbb{R}^1,\,y\ge 0\},$ 
has the function $y(x)\equiv 0$ on $\mathbb{R}^1$ as a solution, according to the provided definitions.}

\smallskip
THe initial value problem with the initial values $x_0,y_0$ will be denoted as IVP\!$(x_0,y_0).$ 

\smallskip
{\bf Df.}\, 
IVP\!$(x_0,y_0)$ for the equation (1) will be called {\it an interior initial value problem} and denoted as IIVP\!$(x_0,y_0),$ 
if a point $(x_0,y_0)\in G,$ and --- {\it a boundary initial value problem} (BIVP\!$(x_0,y_0)),$ if a point $(x_0,y_0)\in \widehat G.$

\smallskip
{\bf Df.}\, 
{\it Interior (boundary, mixed) solution IVP\!$(x_0,y_0)$ exists,}\, if a point $(x_0,y_0)\in G\,(\widehat G,\widetilde G)$ and there is an interval $\langle a,b\,\rangle\ni x_0$ and interior (boundary, mixed) solution $y=\varphi(x),$ defined on this interval, such that $\varphi(x_0)=y_0.$

\smallskip
Thus, the graph of the interior solution of the initial value problem is in the connected open set $G,$ the graph of the boundary solution --- is in $\widehat G,$
the graph of the mixed solution --- is in both set $G$ and set $\widehat G.$ 

Let us present the formulation of the famous Peano theorem on the existence using the provided terminology of the interior initial value problem.

\smallskip
{\bf Peano theorem.} \ {\it Assume the function $f(x,y)$ from the equation (1) is continuous in a connected open set $G\subset \mathbb{R}^2,$ then, at least one solution of the IIVP \!$(x_0,y_0)$ exists for any point $(x_0,y_0)\in G$ and for any Peano segment $P_h(x_0,y_0).$ }

\smallskip
{\small Denote the Peano segment as $P_h(x_0,y_0)=[x_0-h,x_0+h],$ and the constant $h>0$ is defined as follows: 
The constants $a,b>0$ exist, such that the rectangular $R=\{(x,y)\,:\,|x-x_0|\le a,\,|y-y_0|\le b\}\subset G.$
If $f(x,y)\equiv 0$ on $\overline R,$ then $h=a.$
Otherwise, assume \,$M=\max_{(x,y)\in R}|f(x,y)|>0,$ then $h=\min\,\{a,b/M\}.$
Geometrically $h$ is a length of the altitude of the triangle $T^+$ with a vertex at the point $(x_0,y_0),$ which legs have
a slope equal to $\pm M,$ and base belongs to the line $x=x_0+h$ $(T^-$ can be constructed analogously). }

\smallskip
{\bf Remark 3.}\, {\sl It is assumed in paper [7] that the function $f(x,y)$ 
is defined and continuous in the rectangular $R^+=\{(x,y)\,:\,x_0\le x\le x_0+a,\,|y-y_0|\le b\},$ 
which is, in fact, the only requirement for any variant of proving the Peano theorem for $x\ge x_0.$ }

\smallskip
To solve the problem on existence or non-existence of the solution of the boundary initial value problem, we consider, that the equation (1) has form
$$y'=f_0(x,y),\eqno (1_0)$$
where $f_0$ is defined and continuous on the set $\widetilde G=G\cup \widehat G,$ 
where $G$~--- is a connected open set in~$\mathbb{R}^2,$ $\widehat G\subset \partial G,$ 
the point $\text{O}=(0,0)\in \widehat G,\ f_0(0,0)=0,$ and the initial value problem is set at i.\,v.\,$0,0.$

Assume that the initial value problem for the equation (1) is set at the point $(x_0,y_0)\in \widehat G,$ then the change $x=u+x_0,\ \ y=v+y_0+f(x_0,y_0)(x-x_0)$ transforms the equation (1.1) into the equation \,$v'=f_0(u,v),$\, where $f_0(u,v)=f(u+x_0,v+f(x_0,y_0)u+y_0)-f(x_0,y_0).$ For $x=x_0,\ y=y_0,$ we obtain: $u=u_0=0,\ v=v_0=0$ and $f_0(0,0)=0.$ 

The continuity  of the function $f_0(x,y)$ at the boundary point $\text{O}$ means that
$$\forall\,\tau>0\ \ \exists\,\delta_\tau>0\!:\ \forall\,(x,y)\in \overline V_{\delta_\tau}\cap \widetilde G\ \Rightarrow\ |f_0(x,y)|\le \tau\ \ 
  (\overline V_{\delta_\tau}=\{(x,y)\!:\,|x|\le \delta_\tau,\,|y|\le \delta_\tau\}).\eqno (2)$$

{\bf Df.}\, Denote function $y=b_{a_u}^+(x),$ defined on $[0,a_u],$ as {\it a right upper-boundary function}, 
if the five conditions hold:
1) $b_{a_u}^+(x)\in C^1([0,a_u]);$
2)~$b_{a_u}^+(0)=0;$
3) ${b_{a_u}^+}'(0)\ge 0;$
4) $b_{a_u}^+$ is convex on $[0,a_u],$ if ${b_{a_u}^+}'(0)=0;$
5)~{\it the right upper-boundary} curve $\gamma_{a_u}^+=\{x\in [0,a_u],$ $y=b_{a_u}^+(x)\}\subset \widehat G.$  

{\small Here, the symbol \,$+$\, denotes $right,$ $b$~--- $boundary,$ $u$~--- $upper.$ } 

The right lower-boundary curve and its left analogues can be defined in the same manner.

Introduce $\tau_u={b_{a_u}^+}'(0)/2,\ \tau_l=-{b_{a_l}^+}'(0)/2$ and assume without loss of generality, restricting the domain of the boundary functions if needed, that
$$\begin{matrix} b_{a_u}^+(a_u)\le a_u\ \ (\tau_u=0),\quad \forall\,x\in [0,a_u]\!:\ {b_{a_u}^+}'(x)\ge \tau_u\ \ (\tau_u>0); \\
-b_{a_l}^+(a_l)\le a_l\ \ (\tau_l=0),\quad \forall\,x\in [0,a_l]\!:\ \;-{b_{a_l}^+}'(x)\ge \tau_l\ \ \,(\tau_l>0). \end{matrix}\eqno (3)$$

For any chosen $c^*>0$ and $c\ (0<c\le c^*)$ consider the following sets

\smallskip\noindent
$N_c^+=\{(x,y)\!:\,x\in(0,c],\,|y|\le c\},\ O_c^+=\{(x,y)\!:\, x\in (0,a_l],\,b_{a_l}^+(x)\le y\le c;\ x\in (a_l,c],\,|y|\le c\},$

\smallskip\noindent
$U_c^+=\{(x,y)\!:\, x\in (0,a_u],\,-c\le y\le b_{a_u}^+(x);\ x\in (a_u,c],\,|y|\le c\},\ \ B_c^+=U_c^+\cap O_c^+.$
\nopagebreak 

{\small Here, $N$\,means $neighborhood,$\,$U$~---\,$under,$\,$O$~---\,$over,$\,$B$~---\,$between.$ }

\smallskip
Denote that the following cases are realized for the equation $(1_0)$ 

$N^+],$ if $\exists\, c_w>0\!:\ N_{c_w}^+\cap \widehat G=\empty\ \ (w$ ---\,$without);$

$U^+],$ if $\exists\, c_u\ (0<c_u \le c_u^*)\!:\ U_{c_u}^+\cap \widehat G=\gamma_{a_u}^+\backslash \text{O};$ \\
$U_1^+]\!:\ (U_{c_u}^+\!\backslash \gamma_{a_u}^+)\subset G,$ \ two subcases: 
$U_1^{+,>}]\!:\ {b_{a_u}^+}'(0)>0,\ U_1^{+,=}]\!:\ {b_{a_u}^+}'(0)=0;$ \\
$U_2^+]\!:\ U_{c_u}^+\cap G=\emptyset,$ the same subcases; 

$O^+],$ if $\exists\, c_l\ (0<c_l \le c_l^*)\!:\ O_{c_l}^+\cap \widehat G=\gamma_{a_l}^+\backslash \text{O};$ \\
$O_1^+]\!:\ (O_{c_l}^+\!\backslash \gamma_{a_l}^+)\subset G,$ \ two subcases: 
$O_{1,<}^+]\!:\ {b_{a_l}^+}'(0)<0,\ O_{1,=}^+]\!:\ {b_{a_l}^+}'(0)=0;$ \\
$O_2^+]\!:\ O_{c_l}^+\cap G=\emptyset,$ the same subcases;  

$B^+],$ if $\exists\, c_b\ (0<c_b \le c_b^*)\!:\ B_{c_b}^+\cap \widehat G=(\gamma_{a_u}^+\cup \gamma_{a_l}^+)\backslash \text{O};$ \\
$B_1^+]\!:\ (B_{c_b}^+\backslash (\gamma_{a_u}^+ \cup \gamma_{a_l}^+))\subset G,$ \ 4 subcases: 
$B_{1,<}^{+,>}]\!:\ {b_{a_u}^+}'(0)>0,\,{b_{a_l}^+}'(0)<0,$ \\
$B_{1,=}^{+,=}]\!:\ {b_{a_u}^+}'(0),\,{b_{a_l}^+}'(0)=0,$ 
$B_{1,=}^{+,>}]\!:\ {b_{a_u}^+}'(0)>0,\,{b_{a_l}^+}'(0)=0,$ $B_{1,<}^{+,=}]\!:\ {b_{a_u}^+}'(0)=0,\,{b_{a_l}^+}'(0)<0;$ \\
$B_2^+]\!:\ B_{c_b}^+\cap G=\emptyset,$ the same subcases.
	
\begin{center}
\begin{tikzpicture}[scale=0.9]
\tikzset{line01/.style={line width=0.7pt}}
\useasboundingbox (-0.4,-2) rectangle (2.9,2.6);
\fill[whitesmoke] (0,0) ..controls +(0:0.5cm) and +(250:0.6cm).. (2,1.6) -- (2,-2) -- (0,-2) -- cycle;
\draw[->] (-0.3,0) -- (2.8,0) node[below] {\footnotesize $x$};
\draw[->] (0,-2.2) -- (0,2.4) node[left] {\footnotesize $y$};
\draw (0,-2) rectangle (2,2);
\filldraw[line01, draw=black, fill=whitesmoke2] (0,0) ..controls +(0:0.5cm) and +(250:0.6cm).. (2,1.6) node[left=-2pt]
  {\small $\gamma_{\tilde a_u}^+$} -- (2,-2) node[pos=0.28, left=-5pt] {\footnotesize $\overline T_b^+$} -- cycle;
\draw (2.0,-1.2) node {\footnotesize $y=-x$}; 	
\draw[line01, dashed] (0,0) -- (2,0) node[pos=0.5, below=-3pt] {\footnotesize $h_{\vartriangle}^+$};
\draw[dashed] (2,-2) -- (2.2,-2.2);
\draw (0,0)+(-135:7pt) node {\footnotesize 0};
\draw (2.0,-0.25) node {\footnotesize $\tilde a_u=\tilde c$};
\draw (0.1,2) node [left] {\footnotesize $\tilde c$};
\draw (0.1,-2) node [left] {\footnotesize $-\tilde c$};
\draw (1.1,2.3) node {\small $U_1^{+,=}]$};
\end{tikzpicture} \qquad\quad
\begin{tikzpicture}[scale=0.9]
\tikzset{line01/.style={line width=0.7pt}}
\tikzset{line03/.style={line width=1pt}}
\useasboundingbox (0.7,-2) rectangle (2.9,2.6);
\fill[whitesmoke] 
   (0,0) ..controls +(60:0.5cm) and +(210:0.6cm).. (1.75,2) -- (2,2) -- (2,-2) -- (1.65,-2) ..controls +(110:1.7cm) and +(-40:1cm).. cycle;
\draw [->] (-0.3,0) -- (2.8,0) node[below] {\footnotesize $x$};
\draw [->] (0,-2.2) -- (0,2.4) node[left] {\footnotesize $y$};
\draw (0,-2) rectangle (2,2);
\filldraw[line01, draw=black, fill=whitesmoke2] 
   (0,0) -- (intersection of 0,0--20:1cm and 1.75,0--1.75,1) coordinate (t1) -- (intersection of 0,0-- -20:1cm and 1.75,0--1.75,1) coordinate (t2)
	 node[pos=0.3, left=-5pt] {\footnotesize $\overline T_b^+$} -- cycle;
\draw[line01, dashed] (0,0) -- (1.75,0) node[pos=0.75, below=-3pt] {\footnotesize $h_{\vartriangle}^+$};
\draw[dashed] (t1) -- (1.75,2);
\draw[dashed] (t2) -- (1.75,-2);
\draw[dashed] (0,0) -- (30:3cm);
\draw (2.7,1.55) node {\footnotesize $y=\tau_u x$}; %
\draw[line01] (0,0) ..controls +(60:0.5cm) and +(210:0.6cm).. (1.75,2) node[near end,below=+1pt] {\small $\gamma_{\tilde a_u}^+$};
\draw[line01] (0,0) ..controls +(-40:1cm) and +(110:1.7cm).. (1.65,-2);
\draw (1.0,-1.2) node {\small $\gamma_{\tilde a_l}^+$}; %
\draw[dashed] (0,0) -- (60:2.3cm);
\draw (0.2,1.65) node {\footnotesize $y=2\tau_u x$}; 
\draw[dashed] (0,0) -- (-40:3.1cm);
\draw (2.5,-1.3) node {\footnotesize $y=\!-2\tau_l x$}; 
\draw[shift={(t1)},dashed] (0,0) -- (intersection of 0,0--20:1cm and 0.9,0--0.9,1) node[right=2pt, pos=0.1] {\footnotesize $y=\tilde \tau x$};
\draw[shift={(t2)},dashed] (0,0) -- (intersection of 0,0-- -20:1cm and 0.9,0--0.9,1);
\draw (2.7,-0.6) node {\footnotesize $y=\!-\tau_l x$}; 
\draw (0,0)+(-135:7pt) node {\footnotesize 0};
\draw (2,0)+(-45:7pt) node {\footnotesize $\tilde c$};
\draw (0.1,2) node [left] {\footnotesize $\tilde c$};
\draw (0.1,-2) node [left] {\footnotesize $-\tilde c$};
\draw (1.7,0)+(-45:7pt) node {\footnotesize $\tilde a$};
\draw (1.1,2.3) node {\small $B_{1,<}^{+,>}]$};
\end{tikzpicture} \qquad
\begin{tikzpicture}[scale=0.9]
\tikzset{line01/.style={line width =0.7pt}}
\tikzset{line03/.style={line width =1pt}}
\useasboundingbox (-0.4,-2) rectangle (2.9,2.6);
\draw[->] (-0.3,0) -- (2.8,0) node[below] {\footnotesize $x$};
\draw[->] (0,-2.2) -- (0,2.4) node[left] {\footnotesize $y$};
\draw (0,-2) rectangle (2,2);
\filldraw[line01, draw=black, fill=whitesmoke2] (0,0) ..controls +(0:0.5cm) and +(250:0.6cm).. (2,1.6) 
   node[near end, left=0pt] {\small $\gamma_{a_u}^+$} -- (2,-1.7)
	 node[pos=0.3, left=-4pt] {\footnotesize $\overline T_b^+$} ..controls +(100:1cm) and +(0:1.5cm)..
	 node[near start, left=-3pt] {\small $\gamma_{a_l}^+$} cycle;
\draw[line01, dashed] (0,0) -- (2,0) node[pos=0.8, below=-2.5pt] {\footnotesize $h_{\vartriangle}^+$};
\draw (0,0)+(-135:7pt) node {\footnotesize 0};
\draw (1.75,0.15) node {\footnotesize $a_u=a_l=c_b$};
\draw (0.1,2) node [left] {\footnotesize $c_b$};
\draw (0.1,-2) node [left] {\footnotesize $-c_b$};
\draw (1.0,2.3) node {\small $B_{1,=}^{+,=}]$};
\end{tikzpicture} \ \ 
\begin{tikzpicture}[scale=0.9]
\tikzset{line01/.style={line width =0.7pt}}
\tikzset{line03/.style={line width =1pt}}
\useasboundingbox (-0.4,-2) rectangle (1.7,2.6);
\fill[whitesmoke] (0,0) ..controls +(-60:1cm) and +(140:1cm).. (2,-1.8) --
                  (2,1.6) ..controls +(-120:0.5cm) and +(0:0.5cm).. cycle;
\draw[->] (-0.3,0) -- (2.8,0) node[below] {\footnotesize $x$};
\draw[->] (0,-2.2) -- (0,2.4) node[left] {\footnotesize $y$};
\draw (0,2) rectangle (2,-2);
\filldraw[line01, draw=black, fill=whitesmoke2] (0,0) --
   (intersection of 0,0-- -30:1cm and 2,0--2,1) coordinate (t1) --
	 node[pos=0.7, left=-4pt] {\footnotesize $\overline T_b^+$} (2,1.6) ..controls +(-120:0.5cm) and +(0:0.5cm)..
	 node[pos=0.4,left=2pt] {\small $\gamma_{\tilde a_u}^+$} cycle;
\draw[line01, dashed] (0,0) -- (2,0) node[pos=0.6, below=-2.5pt] {\footnotesize $h_{\vartriangle}^+$};
\draw[line01] (0,0) ..controls +(-60:1cm) and +(140:1cm).. (2,-1.8);
\draw (0.4,-1.1) node {\small $\gamma_{\tilde a_l}^+$};
\draw[dashed] (0,0) -- (-60:2.5cm);
\draw (0.95,-1.75) node {\footnotesize $y=-2\tau_l x$}; %
\draw[shift={(t1)},dashed] (0,0) -- (intersection of 0,0-- -30:1cm and 0.6,0--0.6,1);
\draw (2.12,-0.7) node {\footnotesize $y=\!-\tau_l x$};
\draw (0,0)+(-135:7pt) node {\footnotesize 0};
\draw (1.7,0.2) node {\footnotesize $\tilde a_u=\tilde a_l=\tilde c$};
\draw (0.1,2) node [left] {\footnotesize $\tilde c$};
\draw (0.1,-2) node [left] {\footnotesize $-\tilde c$};
\draw (1.1,2.3) node {\small $B_{1,<}^{+,=}]$};
\end{tikzpicture}
\end{center}

The figures illustrate not only the respective sets but also the right boundary triangles and the Peano segments $(\tau_u,\tau_l$ from (3)).
The figures illustrating the other five cases can be found in [1,\,\S\,3]. 

Consider the situation, that may arise, for example, in cases $U_1^{+,=}].$

For any small $\tau$ in the formula (2) and $T_b^+$ in $N_{\delta}^+,$ there is always the neighborhood of the point ${\text O},$ such that the right upper-boundary curve $\gamma_{c_u}^+\ (c_u=\delta)$ is located under the upper leg of the triangle.
Therefore, the Euler polygon can't be continued right, if it reaches $\gamma_{c_u}^+$ at some point $(x_*,b_{c_u}^+(x_*)),$ if the segment of the slope field will have the slope lesser or equal to $\tau,$ but greater than the slope of the tangent to the $\gamma_{c_u}^+$ at this point, which is possible (the polygon can't leave $\widetilde G).$ 
\begin{center} \includegraphics[scale=0.6]{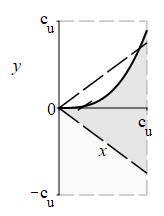} \end{center}  

To mend this problem, introduce restrictions on the functions $f_0$ at any point of the curves $\gamma_{a_u}^+,\,\gamma_{a_l}^+$ in cases $U_1^{+,=}],\,O_{1,=}^+],\,B_{1,=}^{+,=}],\,B_{1,=}^{+,>}]$ ш $B_{1,<}^{+,=}]:$ 

$$\begin{matrix}
  \phantom{aaai} \forall\,x\in (0,a_u]\!:\ f_0(x,b_{a_u}^+(x))\le {b_{a_u}^+}'(x), \text{ хёыш } {b_{a_u}^+}'(0)=0; \\
  \phantom{aaaa} \forall\,x\in (0,a_l]\!:\ \ f_0(x,b_{a_l}^+(x))\ge {b_{a_l}^+}'(x),\ \text{ хёыш }\ {b_{a_l}^+}'(0)=0,\end{matrix}\eqno(4)$$
which mean, that the right half-segment of the slope field of the equation $(1_0)$ is directed inside or along the boundary of a connected open set $G$ at any point of $\gamma_{a_u}^+$ and $\gamma_{a_l}^+.$

\smallskip
{\bf Theorem}\, (on existence of the solution of the boundary initial value problem).\, {\it
Assume that in equation~$(1_0)$ the function $f_0(x,y)$ is continuous on the set $\widetilde G,$
then in cases $U_1^{+,=}],\,O_{1,=}^+],\,B_{1,=}^{+,=}],\,B_{1,=}^{+,>}],\,B_{1,<}^{+,=}]$ at least one solution of BIVP\!$(0,0)$ exists on any right boun- dary Peano segment if the respective condition $(4)$ holds.}

Prove can be found in [1,\,p.\,4]. The theorem on non-existence of the solutions of the BIVP is also proven there. 

\smallskip
The questions, related to the extendability of the solutions of the boundary initial value problem, are provided in the academic book [6,\,ch.\,1,\,\S\,1,\,p.$5^0,6^0$].

\smallskip
\begin{center} {\large 2.\, On uniqueness or non-uniqueness of the solution \\ of the boundary initial value problem.} \end{center}  
\noindent\!2.1. {\bf Non-uniqueness points.} 

${\bf Df.}$\, Denote point $(x_0,y_0)\in \widetilde G$ as {\it a non-uniqueness point,} 
if the solutions $y=\varphi_1(x)$ and $y=\varphi_2(x)$ of the IVP\!$(x_0,y_0),$ defined on $\langle a,b\,\rangle \ni x_0,$ and a sequence $x_k\to x_0$ when $k\to \infty,$ $x_k\in \langle a,b\,\rangle $ exist, such that $\varphi_1(x_k)\ne \varphi_2(x_k)$\ $(k=1,2,\ldots).$
Otherwise, $(x_0,y_0)$ is a {\it formal uniqueness point.}

\smallskip
Let us provide the definition of the non-uniqueness point, equivalent to the first one. Its negation allows to easily formulate the <<straight>> definition of the formal uniqueness point, which is:\,
{\sl denote point $(x_0,y_0)\in \widetilde G$ as a {\it non-uniqueness point,}
if such solutions  $y=\varphi_1(x)$ and $y=\varphi_2(x)$ of the IVP\!$(x_0,y_0),$ defined on $\langle a,b\,\rangle \ni x_0,$ exist, that, for any interval $(\alpha,\beta)\ni x_0,$ such point $x^*\in (\alpha,\beta)\cap \langle a,b\,\rangle $ exists, that $\varphi_1(x^*)\ne \varphi_2(x^*).$ }

\smallskip
It is important here, that $(\alpha,\beta)$ is an interval, because it eliminates the case, when $(\alpha,\beta\,\rangle =(\alpha,x_0],$ but $\langle a,b\,\rangle =[x_0,b\,\rangle ,$ and the point $x^*$ can't exist.

\smallskip 
Using common reasoning, it may be stated that the integral curves <<branch>> at the non-uniqueness point, 
but the geometrical understanding of the term {\it branching} show that it's a~wrong assumption. 

\smallskip
{\bf Df.}\ The integral curve {\it branches from the right} at the point $(x_0,y_0)\in \widetilde G,$
if such solutions $y=\varphi_1(x)$ $y=\varphi_2(x)$ of the IVP\!$(x_0,y_0)$ and the constant $\delta>0$ exist that for any $x\in (x_0,x_0+\delta),$  
$\varphi_1(x)\ne \varphi_2(x)$\, (branching from the left can be defined analogously).

\smallskip
{\bf Example 1.}\, 
Consider the equation \,$y'=3\sqrt{x}\sqrt{y}.$\,  
Its right-hand side is continuous on $\widetilde G=\{(x,y)\!:\,x\ge 0,\ y\ge 0\},$ 
the set $\widehat G$ consists of rays $y\equiv 0$ $(x\ge 0)$ and $x\equiv 0$ $(y>0).$

The function $y(x)\equiv 0,$ $x\in [0,+\infty)$ is a boundary solution, all interior solutions 
can be defined, using formula $\varphi(x)=(x^{3/2}-C)^2,$ where $x>0$ when $C\le 0$ and $x>C^{2/3}$ when $C>0.$ 
Its graphs are in the connected open set $G=\{(x,y)\!:x>0,\,y>0\}.$ 

All interior solutions with $C\le 0$ become mixed solutions, if the point  $x=0$ is added to their domain, i.\,e. for any constant $C\le 0,$ the function $y=(x^{3/2}-C)^2$ on $[0,+\infty)$ is a mixed solution and $(0,C^2)$ is its only boundary point.
Moreover, for any constant $C>0,$ 
the <<piece-wise>> function $y=\{0\ \hbox{ яЁш }\ 0\le x\le C^{2/3},$ $(x^{3/2}-C)^2\ \hbox{ when }\ x>C^{2/3}\}$ is a mixed solution.

Therefore, the boundary solution $y\equiv 0,$ set, for example, on $[0,1]$ can be continued to $\widetilde b,$ for any $\widetilde b>1,$ 
while keeping being boundary. At the same time, for $x_*\in [1,\widetilde b),$ this solution can become 
mixed solution $y(x)=\{\,0\ \hbox{яЁш}\ x\in [0,x_*^{2/3}],$ $(x^{3/2}-x_*)^2\ \hbox{when}\ x\in (x_*^{2/3},\widetilde b\,\rangle \}.$ \ $\square$

\includegraphics[scale=0.6]{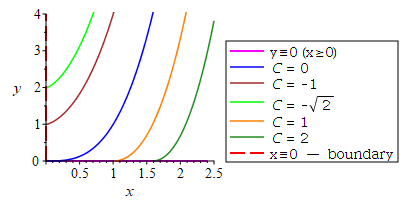}   

Example 1 demonstrates a typical situation, when the integral curves branch.
For any $x_*\ge 0,$ the integral curve of the solution $y=(x^{3/2}-x_*^{3/2})^2,\ x\ge x_*$ branches from the right from the graph of the boundary solution at the point $(x_*,0)$ and doesn't have any common points with the abscissa axis.

{\bf Counterexample 1.}\, 
Introduce some notations: 
$$\begin{matrix} 
  \hfill a_n=2^{-n},\ \ b_n=2^{-2(n+1)},\ \ d_n=3\cdot 2^{-(n+1)},\ \ J_n=[2^{-n},2^{-(n-1)}]; \\ 
	\hfill \psi_n(x)=\big(b_n-(x-d_n)^2 \big)^{3/2},\ x\in J_n;\ \ \ D_n=\left\{(x,y)\!:\,x\in J_n,\,|y|\le \psi_n(x)\right\}, \\ 
	  \hfill D_n^+=\left\{(x,y)\!:\,x\in J_n,\,y\ge \psi_n(x)\},\ \ D_n^-=\{(x,y)\!:\,x\in J_n,\,y\le -\psi_n(x)\right\}; \\  
  h_n(x,y)=-3(x-d_n)\cdot \left\{ y^{1/3}\hbox{ when } (x,y)\in D_n,\ \pm \psi_n^{1/3}(x)\hbox{ when } (x,y)\in D_n^{\pm}\right\}  
  \end{matrix} \quad\ (n\in \mathbb{N});$$
$$\begin{matrix}
  \hfill \psi(x)=\big\{ \psi_1(x),\,x\in J_1;\ \psi_2(x),\,x\in J_2;\,\ldots \big\},\ \ \psi(0)=0; \\ 
	h(x,y)=\big\{ h_1(x,y),\,x\in J_1;\ h_2(x,y),\,x\in J_2;\,\ldots \big\},\ \ h(0,y)\equiv 0.\end{matrix} $$

Let us note, that $d_n$ is a middle-point of segment $J_n,$ $b_n\ge (x-d_n)^2$ when $x\in J_n;$
$\psi_n(a_n)=\psi_n(a_{n-1})=0,$ $\max_{x\in J_n} \psi_n(x)=\psi_n(d_n)=b_n^{3/2}=2^{-3(n+1)}\to 0$ when $n\to \infty;$
the functions$\psi(x),\,\psi'(x),\,h(x,y)$  are continuous when $x\in [0,1],$ $\psi'(0)=0.$

Consider the equation 
$$y'=h(x,y),\eqno (5)$$
defined and continuous on the set $\widetilde G=\{(x,y)\!:\,x\ge 0,\,y\in \mathbb{R}^1\}$ 
and invariant with respect to change $y$ to $-\tilde y,$ therefore, it's sufficient to solve it when $y\ge 0.$ \goodbreak

The equation (5) takes the following form on the set $D_n$
$$y'=-3(x-d_n)y^{1/3}.$$

By integrating, we obtain the trivial solution $y(x)=0,\ x\in J_n$ and the general integral $U_n(x,y)=C$ with $U_n=y^{2/3}+(x-d_n)^2,$ 
where $C\in (0,b_n].$ The integral $U_n(x,y)=b_n,$ while being solved with respect to $y,$ 
coincides with the function $y=\psi_n(t),\ t\in J_n,$ and the function $U_n(x,y)=0$ degenerates into the point $(d_n,0)$ when $C=0.$ 

The equation (5) takes the following form on the set $D_n^+$
$$y'=-3(x-d_n)\psi_n^{1/3}(x).$$

Because $\psi_n'(x)\equiv -3(x-d_n)\psi_n^{1/3}(x)\,(\!=-3(x-d_n)(b_n-(x-d_n)^2 )^{1/2}$ and $(\psi_n(x)+c)'=\psi_n'(x),$ 
the function $\psi_n(x)+c\ \,(c\ge 0)$ is its general solution. 

Moreover, because $\psi_n''(x)=0\Leftrightarrow b_n-2(x-d_n)^2=0,$ the function $\psi_n'(x)$ reaches its extreme values at the points $x_n^{\pm}=d_n\pm (b_n/2)^{1/2}=(3\pm 2^{-1/2})2^{-(n+1)}\in J_n.$ Because $x_n^{\pm}-d_n=\pm 2^{-n-3/2}$ and $\big(b_n-(x-d_n)^2 \big)^{1/2}=2^{-n-3/2},$ $\psi_n'(x_n^{\pm})=\mp 3\cdot 2^{-2n-3},$ which means that the maximal absolute value of the ratio of the slope of the tangent of the solution's $y=\psi_n(x)$ integral curve to the length of the segment $J_n$ is equal to $(3/8)\cdot 2^{-n}.$

As a result, the equation (5) has a trivial solution $y(x)=0,$ two sets of solutions: $\varphi^{\pm}(x,c)=\pm(\varphi(x)+c)$ $(c\ge 0),$ defined on the segment $[0,1]$ and, for any $n\in \mathbb{N},$ two sets of solutions: $\varphi^{\pm}(x,C)=\pm\big(C-(x-d_n)^{2/3}\big)^{3/2}\ (C\in (0,b_n)),$ obtained from the general integral $U_n(x,y)=C.$ Graphs of these two sets are in ${\rm Int}\,D_n\backslash \{(d_n,0)\}.$ \ $\square$

\includegraphics[scale=0.42]{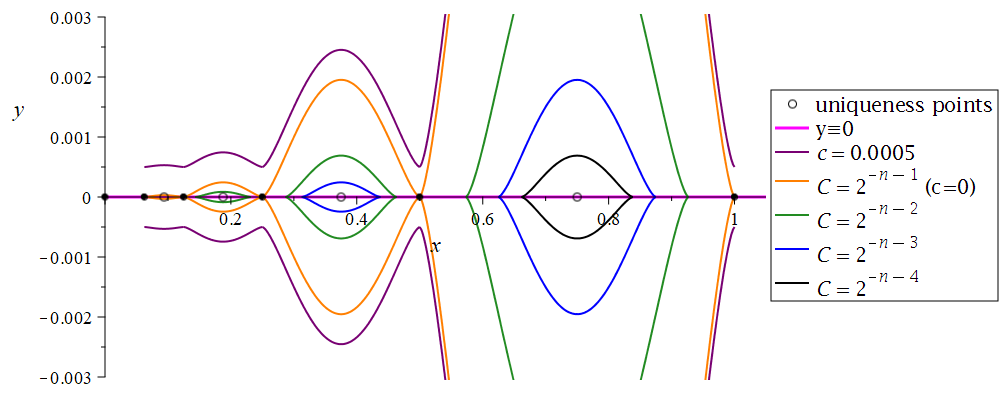}

\smallskip
Thus, the counterexample \;1 demonstrates, that the point $(0,0)$ is a non-uniqueness point by definition, but it can't be stated, that the integral curve of the solutions $y=\pm \varphi(x)$ of the IVP\!$(0,0)$ branch from the integral curve of the solution $y(x)=0.$ 

\smallskip
{\bf Remark 4.} \ 
{\sl Function $h(x,y)$ in the equation (5) can be easily redefined in such way, that is becomes continuous in $\mathbb {R}^2$ and invariant with respect to change $x$ to $-\tilde x.$ It means, that a non-uniqueness point, which is not a branching point, can be both boundary and interior. }

\smallskip
\noindent\!2.2. {\bf Formal and local uniqueness of the boundary initial value problem.}  
The negation of the definition of the <<non-uniqueness point>> is different for boundary and interior points from $\widetilde G,$ because, unlike interior case, the solution of the BIVP\!$(x_0,y_0)$ may not exist or may be not extendable left or right beyond the point $x_0.$

\smallskip
{\bf Df.}\ {\it Denote point $(x_0,y_0)\in \widehat G$} as {\it a formal uniqueness point,} 
and {\it solution of the BIVP\!$(x_0,y_0)$} of the equation (1) as {\it a formally unique (or unique at the point),} 
If this initial value problem has no solution or, for any two its solutions $y=\varphi_1(x)$ and $y=\varphi_2(x),$ defined on some interval $\langle a,b\,\rangle ,$ there is an interval $(\alpha,\beta)\ni x_0,$ such that, for any $x\in (\alpha,\beta)\cap \langle a,b\,\rangle ,$ $\varphi_1(x)=\varphi_2(x).$

\smallskip
The use of the word <<formal>> in the definition of the solution of the BIVP\!$(x_0,y_0)$ is not accidental,
because there is no information about how close the abscissa $x_1$ of the closest to $(x_0,y_0)$ non-uniqueness point $(x_1,\varphi(x_1))$ to the $x_0.$

Indeed, by extracting some specific solution of the initial value problem, set at the formal uniqueness point, we conclude, that, by definition, the other solutions of this initial value problem, if they exist, coincide with the extracted one on its intervals, while the length of these intervals may converge to zero.

It is exactly the reason, why suggested definition of uniqueness is not sufficient from a standpoint of the applied mathematics, because the process, described by the solution $y=\varphi(x)$ BIVP\!$(x_0,y_0),$ is not deterministic, for example, in such case, when, for any $x>x_0,$ the arc of the integral curve of another solution branches from the graph of the solution $y=\varphi(x)$, while the point $(x_0,y_0)$ is the only non-branching point.

Counterexample 2 and example 2 will show, that the described situation takes place, i.\,e. the point $(x_0,y_0)$  is a formal uniqueness point due to the natural restriction on the set $\widetilde G,$ which is the domain of the right-hand side of the equation (1). When the function $f(x,y)$ is continuously redefined to some extension of the set $\widetilde G,$ the point $(x_0,y_0)$ becomes the non-uniqueness point, which compels to denote it as <<hidden>> non-uniqueness point.
But, as it will be stated further, such situation may take place only for boundary point of the set $\widetilde G.$

Therefore, the different definition of uniqueness of the solution of the BIVP\!$(x_0,y_0)$ is more informative. It also guarantees that the process on some common segment, containing the point $x_0,$ is deterministic.  

\smallskip
{\bf Df.}\, {\it Denote solution of BIVP\!$(x_0,y_0)$} of the equation (1) 
as {\it a right locally unique,} if it exists on a segment $[x_0,b),$ and such segment $[x_0,\beta)$ exists, that all solutions of this initial value problem are extendable to $[x_0,\beta),$ and, for any two solutions $y=\varphi_1(x)$ and $y=\varphi_2(x),$ extended, if needed, to $[x_0,\beta),$ the following takes place: $\varphi_1(x)\equiv \varphi_2(x)$ on $[x_0,\beta).$ Otherwise, the solution of BIVP is not a right locally unique.

{\it Left locally unique solution} is defined analogously.

\smallskip
{\bf Df.}\, {\it Denote solution of the BIVP\!$(x_0,y_0)$} as {\it locally unique,} and point $(x_0,y_0)$ as {\it a uniqueness point,} 
if this solution is left locally unique and can't be extended to the right beyond the point $x_0,$ or it is right locally unique 
and can't be extended to the left beyond the point $x_0,$ or both left and right locally unique.
Otherwise, the solution of BIVP \!$(x_0,y_0)$ is not locally unique, and the point $(x_0,y_0)$ is not a uniqueness point.

\smallskip
In example\;1, for any point $y^*>0,$ the solution of BIVP\!$(0,y_*)$ takes form $y=\big(x^{3/2}+y_*^{1/2}\big)^2,$ defined on interval $[0,+\infty),$ and is right locally unique, and, for $y_*=0,$ the boundary solution $y(x)=0$ branches from the solution $y=x^3$ of the BIVP\!$(0,0),$ which means that the solution is not right locally unique and that origin is non-uniqueness point.
In counterexample\;1, there is no $\partial h(x,y)/\partial y$ in the equation (5) when $y=0,$ 
but the graph of the solution $y(x)=0$ contains the countable number of uniqueness points $(d_1,0),(d_2,0),\ldots\,.$ 

\smallskip
Next example demonstrates that the point $(x_0,y_0)\in \widehat G,$ while being the formal uniqueness point, may be not uniqueness point, i.\,e. 
for any two solutions of BIVP\!$(x_0,y_0),$ there is an interval $[x_0,b),$ on which the graphs of these solution coincide, but there is no interval, where the graphs of all solutions coincide. Therefore, there is no right local uniqueness of the solution of BIVP at such point.

\smallskip
{\bf Counterexample\;2.} \ 
Consider homogeneous equation
$$y'=\sqrt{y}-2\sqrt{x^2-y}+x\eqno (6)$$
defined on $\widetilde G=\{(x,y)\!:\, x\ge 0,\ 0\le y\le x^2\},$ then its boundary is $\widehat G=\gamma_1\cup\gamma_2,$
where $\gamma_1=\{x\ge 0,\ y=x^2\},\ \gamma_2=\{x>0,\ y=0\}.$ The function $y=x^2,$ when $x\in [0,+\infty)$ is a boundary solution. 

The change $y=u^2x^2,$ where $u\in [0,1],$ transforms equation (6) into equation $2xuu'+h(u)=0,$
in which $h=2\sqrt{1-u^2}-2(1-u^2)+(1-u),$ $h(u)>0$ when $u\in [0,1),$ $h(1)=0$ 
and equality $u\equiv 1$ defines the boundary solution $y=x^2,\ x\ge 0.$ 
Using the variable separation method, the solution $\displaystyle \int_0^u \frac{2v}{h(v)}\,dv+\ln\frac{x}{c_*}=0$ is found for any constant $c_*>0.$

Substituting the variable $u\in [0,1)$ in the upper integral limit with the function $x^{-1}\sqrt{y}$  
and writing down the general solution in the form of the general integral, we obtain: 
$$U(x,y)=c\ \ (x>0,\ y\ge 0),$$ 
where $U=xe^{-\eta(x,y)},$ $\displaystyle \eta(x,y)=\int_0^{\sqrt{y}/x} \frac{2v}{h(v)}\,dv\ \ (x>0),$ $h=2\sqrt{1-v^2}-2(1-v^2)+(1-v).$ 

\includegraphics[scale=0.56]{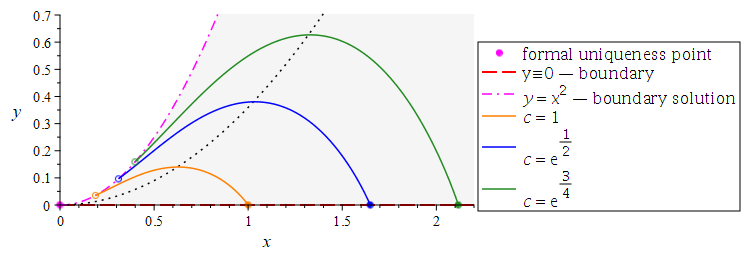}  

It's important to notice, that the integral with a variable upper limit, while having singularity when $v=1,$ is proper, and $\displaystyle \theta=\int_0^1 \frac{2v}{h(v)}\,dv\approx 1.66.$
Therefore, for any $\tilde c>0,$ the arc of the integral curve of the solution $U(x,y)=\tilde c$\, as value of $x$ decreases, comes into contact with the integral curve $\gamma_1,$ when $x_{\tilde c}=\tilde ce^{-\theta},$ instead of converging to the point $(0,0)$ along the curve.
Also, $U(\tilde c,0)=\tilde c$ for any point $(\tilde c,0)$ of the lower bound, because $\eta(x,0)\equiv 0\ \ (x>0),$ and $U(\tilde c,\tilde c^2)=\tilde c e^{-\theta}$ for any point $(\tilde c,\tilde c^2)$ of the upper bound, because $\eta(x,x^2)\equiv \theta.$ \ $\square$

\smallskip
This example demonstrates the following: if, for any $\tilde c>0,$ integral $U(x,y)=\tilde c$ is written down in form $y=\varphi(x,\tilde c),$ 
which is theoretically possible, then the piece-wise function $y=\{x^2,\ x\in [0,x_{\tilde c}];\ \varphi(x,\tilde c),\ x\in [x_{\tilde c},\tilde c]\},$ where $x_{\tilde c}=\tilde ce^{-\theta},\,\varphi(\tilde c,\tilde c)=0,\,\varphi(x_{\tilde c},\tilde c)=x_{\tilde c}^2,$ 
is a complete mixed solution of both BIVP\!$(\tilde c,0),$ and BIVP\!$(0,0).$

In other words, the origin is a formal uniqueness point, while any boundary point $(x_c,x_c^2)$ \ $(x_c>0)$ is a non-uniqueness point, 
because the arc of the integral curve branches from the graph of the boundary solution $y=x^2$ at it, and the arc reaches the lower bound 
at the point $(c,0)$ as the value of $x$ increases. Therefore there is no interval for solutions of BIVP\!$(0,0),$ where they all coincide.

\smallskip
{\bf Example \;2.} \ 
Consider now the <<extended>> equation
$$y'=\{ \sqrt{y}-2\sqrt{x^2-y}+x\,\hbox{ when }\, 0\le y\le x^2,\ \sqrt{x^2+y}-2x\,\hbox{ when }\, -x^2\le y\le 0 \},\eqno (6_e)$$ 
its right-hand side is continuous on $\widetilde G_e=\{(x,y)\!:\, x\ge 0,\ y\in \mathbb{R}^1\}$ 
and coincides with the right-hand side of the equation (6) on $\widetilde G.$

Let us integrate the equation $(6_e)$ at the lower semi-plane and <<conjoin>> the obtained solution with the solution of the equation (6) 
on the abscissa axis.

Assume $y\le 0.$ Then the change $y=u^2-x^2$ $(0\le u\le x)$ transforms the equation $y'=\sqrt{x^2+y}-2x$ into the equation $(2u'-1)u=0.$
Then, $u\equiv 0$ or $2u=x+c,$ therefore,  
$$y=-x^2,\ x\ge 0 \hbox{ --- boundary solution};\quad V(x,y)=c \hbox{ --- general integral},$$
where $V=2\sqrt{x^2+y}-x\quad (x\ge 0,\ y\le 0).$ 

In the first quadrant the solutions of the equation $(6_e)$ are general integral $U(x,y)=c$ and boundary solution $y=x^2,$ obtained in the example\;2.

\smallskip 
\includegraphics[scale=0.54]{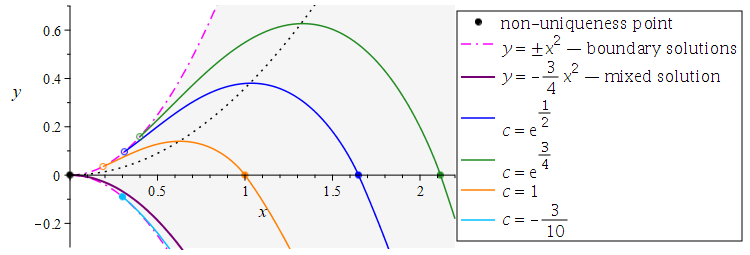}  

Let us note the following:

a) for any $\tilde c>0,$ we have:  $V(\tilde c,0)=\tilde c=U(\tilde c,0),$ 
therefore, the integral curve, parametrized by the integral $U(x,y)=\tilde c$ in the upper semi-plane, 
and the integral curve, defined by the integral $V(x,y)=\tilde c$ in the lower semi-plane, "stitch together" at any point $(\tilde c,0);$

b) integral $V(x,y)=0$ defines the mixed solution $y=-3x^2/4,\ \ x\ge 0;$ 

c) for any constant $\tilde c<0,$ inexplicitly given solution $2\sqrt{x^2+y}=x+\tilde c$ is defined when $x\ge -\tilde c$ 
and its graph is in lower semi-plane;

d) when $x\ge \max\{0,-c\}$ the general integral $V(x,y)=c$ can be presented in a form of the function $y(x,c)=-(3x^2-2cx-c^2)/4,$ 
which has zeroes at $-c/3,\,c.$ \ $\square$

\smallskip 
Thus, example\;2 demonstrates, that the origin is a non-uniqueness point for the equation $(6_e),$ 
because the mixed solution $y=-3x^2/4$ and the boundary solution $y=-x^2/2$ of BIVP\!$(0,0)$ appeared.

Provided examples require some clarification in terminology. 

\smallskip
{\bf Df.}\, Denote formal uniqueness point, which is not uniqueness point, as {\it hidden non-uniqueness point. }

The set of formal uniqueness points consists of the set of uniqueness point and the set (possible empty) of hidden non-uniqueness points.

\medskip 
\noindent\!2.3. {\bf Formal and local uniqueness of the interior initial value problem.} 
In this section we will study the uniqueness of the interior points of the set $\widetilde G,$ taking into account that, for any point $(x_0,y_0)\in G,$ according to the Peano theorem, at least one solution, defined on the Peano segment $P_h(x_0,y_0)=[x_0-h,x_0+h]$ exists. 

\smallskip
{\bf Df.}\, {\it Denote point $(x_0,y_0)\in G$} as {\it a formal uniqueness point,} 
if for any two solution of IIVP\!$(x_0,y_0)$ $y=\varphi_1(x)$ ш $y=\varphi_2(x)$ 
there is such interval $(\alpha,\beta)\ni x_0,$ that  $\varphi_1(x)\equiv \varphi_2(x)$ on $(\alpha,\beta).$ 

\smallskip
{\bf Df.}\, {\it Denote solution of IIVP\!$(x_0,y_0)$} of the equation $(1)$ as {\it unique} or, which is the same, 
{\it locally unique} and point $(x_0,y_0)$ as {\it a uniqueness point,} if such interval $(\alpha,\beta)\ni x_0$ exists, that
all solutions of BIVP\!$(x_0,y_0)$ are extendable on $(\alpha,\beta)$ and, for any two solutions $y=\varphi_1(x)$ and $y=\varphi_2(x),$
extended on $(\alpha,\beta)$ if needed, following condition holds: 
$\varphi_1(x)\equiv \varphi_2(x)$ on $(\alpha,\beta).$
Otherwise, {\it solution of IIVP\!$(x_0,y_0)$ is not unique.} 

\smallskip
Let us show, that, for any two solutions of IIVP\!$(x_0,y_0),$ if such interval $(\alpha,\beta)\ni x_0$ exists, 
that these solution coincide on it (different interval for different pairs), then <<the universal interval>> exists, 
such that all solutions of the initial value problem, extended on this interval if needed, coincide on it.

\smallskip
{\bf Proposition\;1.}\, 
{\it Assume $y=\varphi(x)$ is a solution of IIVP\!$(x_0,y_0),$ 
defined on the Peano segment $P_h(x_0,y_0).$ Then any other solution $y=\psi(x)$ of the same initial value problem,
defined on the interval $\langle a,b\,\rangle \subsetneqq [x_0-h,x_0+h],$ is extendable on $P_h(x_0,y_0).$ }

\smallskip
This proposition allows to choose any interval from $\overline P_h(x_0,y_0)$ containing $x_0$ in the definition of local uniqueness 
of the solution of IIVP\!$(x_0,y_0).$

\smallskip
{\bf Proposition 2.}\, {\it
Assume $\{\chi_k(x)\}_{k=1}^\infty$ is an arbitrary sequence of solutions of IIVP\!$(x_0,y_0),$ defined on some Peano segment $P_h(x_0,y_0).$ Then
$$\forall\,k\in \mathbb{N},\ \forall\,x\in [x_0-h,x_0+h]\!:\ \ 
  \chi_k^l(x)=\min\{\chi_1(x),\ldots,\chi_k(x)\},\ \chi_k^u(x)=\max\{\chi_1(x),\ldots,\chi_k(x)\}$$
is also the solution of the same initial value problem on $P_h(x_0,y_0).$ }

\smallskip
{\bf Lemma}\, (on {\sl lower} and {\sl upper} solutions).
{\it Such solutions of IIVP\!$(x_0,y_0)$ $y=\chi^l(x)$ and $y=\chi^u(x)\ \ (l-lower,\,u-upper)$ exist, that }
$$\forall\,k\in \mathbb{N},\ \forall\,x\in [x_0-h,x_0+h]\!:\ \ \chi^l(x)\le \chi_k^l(x),\ \chi^u(x)\ge \chi_k^u(x).\eqno (7)$$

P r o o f\,. \
Consider, for example, the sequnece of solutions $\{\chi_k^l(x)\}_{k=1}^\infty,$ defined on segment $[x_0,x_0+h].$
This sequence is uniformly bounded and equicontinuous, because the graphs of all these solutions are located in the Peano triangle $T^+.$ Therefore, according to the Arzela-Ascoli theorem, the subsequence exists and it converges uniformly on $P_h(x_0,y_0).$ The limit of this such subsequence is also a solution of the equation (1) on the Peano segment. But sequence $\chi_k^l(x)$ decreases monotonically, 
therefore it converges to the solution $y=\chi^l(x),$ called {\sl lower solution}, for which inequality (7) holds.

The reasoning for the segment $[x_0-h,x_0]$ and the proof of the convergence of the functions $\chi_k^u(x)$ to {\sl upper solution}
$y=\chi^u(x)$ is analogous.\ \ $\blacksquare$

\smallskip
{\bf Theorem}\, (on local uniqueness pf the solution of the interior initial value problem).\,
{\it Assume point $(x_0,y_0)\in G$ is a formal uniqueness point,
then a solution of IIVP\!$(x_0,y_0)$\,of the equation $(1)$ is locally unique. } 

P r o o f\,. \ 
Assume point $(x_0,y_0)\in G$ is a formal uniqueness point. Construct an arbitrary Peano segment $P_h(x_0,y_0).$ 
Reasoning from the contrary, assume that, for any interval  $(\alpha,\beta),$ such that $x_0\in (\alpha,\beta)\subset [x_0-h,x_0+h],$ a solutions of IIVP\!$(x_0,y_0)$ $y=\varphi(x)$ and $y=\psi(x)$ exist and don't coincide on $(\alpha,\beta).$ 

Then, for any $k=1,2,\ldots,$\, the solutions $y=\varphi_k(x)$ and $y=\psi_k(x),$ defined on $P_h(x_0,y_0),$ of this initial value problem exist, such that
$$\exists\,x_k\in (x_0-h/k,x_0+h/k) :\ \ \ \varphi_k(x_k)<\psi_k(x_k).$$

According to proposition\;2, $\varphi_k^l=\min\{\varphi_1(x),\ldots,\varphi_k(x)\}$ and $\psi_k^u=\max\{\psi_1(x),\ldots,\psi_k(x)\}$ are solutions, and according to lemma on {\sl lower} and {\sl upper} solutions, solutions $y=\varphi^l(x)$ and $y=\psi^u(x)$ exist and satisfy the inequation simialr to (7).

As a result, $x_k\to x_0$ when $k\to +\infty$ and the following inequations are valid
$$\forall\,k=1,2,\ldots:\ \ \ \varphi^l(x_k)\le \varphi_k(x_k)<\psi_k(x_k)\le \psi^u(x_k),$$
which means, that $(x_0,y_0)$ is non-uniqueness point. This is a contradiction. \ $\blacksquare$

\smallskip
{\bf Remark 4.}\, {\sl 
Different proof of this theorem is widely known ([7,\,ch.\,III, \S\,2,\,th.\,2.1]).
This proof establishes the existence of so called minimal $y=\underline \chi(x)$ and maximal $y=\overline \chi(x)$ solutions 
of the interior initial problem value on the Peano segment, which possess the property, 
that for any other solution $y=\chi(x)$ of the same initial value problem the following inequalities are valid
$\forall\,x\in [x_0-h,x_0+h]:\ \ \underline \chi(x)\le \chi(x)\le \overline \chi(x).$

As a result, the definition of formal uniqueness of the solution of IIVP\!$(x_0,y_0)$ implies the existence of the interval $(\alpha,\beta)\ni x_0,$ such that the $\underline \chi(x)\equiv \overline \chi(x)$ and, hence, all other solution of the initial value problem coincide on $(\alpha,\beta).$ }

\smallskip
Thus, for all points from a connected open set $G$ the definitions of uniqueness point and formal uniqueness point are equivalent, while the definitions of the local and formal uniqueness are also equivalent for a solution of any IIVP.

\medskip 
\noindent\!2.4. {\bf Theorems on local and global uniqueness of the solutions of BIVP.} 
Provided results explain, why the analogous theorem can't be proven for the BIVP, and that the situation, 
recreated in the counterexample \;2, is possible.

The reason is that to proof the existence of lower and upper solutions (see \; lemma on lower and upper solution), there must be a possibility to construct, for example, right boundary Peano segment, which allows to apply the theorem on existence of a right boundary solution, and this segment, unlike the Peano segment, which exists for any point of a connected open set, may not exist.

In fact, the function $f(x,y)=\sqrt{y}-2\sqrt{x^2-y}+x$ from the equation (5) is continuous on $\widetilde G=\{(x,y)\!:\, x\ge 0,\ 0\le y\le x^2\},$ for which the case $B_{1,=}^{+,=}],$ is realized. Also, the condition (4) doesn't hold for the lower boundary curve  $y(x)=0,$ because $f(x,0)=-x<0,$ and holds for the upper boundary curve $y=x^2,$ because $f(x,x^2)=2x.$ 

Therefore, it is not possible to construct the right boundary Peano triangle. Though the solution of the initial value problem, set at the origin, still exists, the lower solution of the equation (5) does not (upper solution is $y=x^2).$

The method, that was used in the paper [2] in the study of the existence of the BIVP solutions, was used to construct this solution.

Continuous redefinition of the function $f$ in the fourth quadrant removes the boundary $y\equiv 0,$ which <<cuts out>> the lower solution. 
As a result, the <<extended>> equation $(5_e)$ has a lower solution $y=-3x^2/4$ and a boundary solution $y=-x^2,$ 
but the origin becomes non-uniqueness point. This is exactly the reason, 
why the point $(0,0)\in \widehat G$ in equation (5) is a hidden non-uniqueness point.

\smallskip
{\bf Theorem}\, (on local uniqueness of the solution of the boundary initial value problem).\, {\it
Assume the point $(x_0,y_0)\in \widehat G$ is a formal uniqueness point in the equation $(1),$ 
then a solution of BIVP\!$(x_0,y_0)$ is right locally unique, if after the transformation of the equation $(1)$ into the equation $(1_0)$ 
one of the nine cases, considered in the theorem on the existence of the boudnary solution, is realized.}

P r o o f\,.\, 
Assume that any of the cases $U_1^{+,>}],\,O_{1,<}^+],\,B_{1,<}^{+,>}]$ 
or any of the cases $U_1^{+,=}],\,O_{1,=}^+],$ $B_{1,=}^{+,=}],\,B_{1,=}^{+,>}],\,B_{1,<}^{+,=}]$ under condition (4) 
for the BIVP\!$(0, 0)$ of the equation $(1_0)$ is realized. Then, in each case, the right boundary triangles $T_b^+$ 
and Peano segments $P_{h^+}^+(0,0)=[0,h^+]$ $(h^+>0)$ can be constructed for the point $(0,0).$

According to the theorem on existence of the solution of BIVP, at least one solution $y=\varphi(x)$ exists on any such interval, that $\varphi(0)=0.$ 
Local uniqueness of this solution is proven exactly the same way as the theorem on local uniqueness of the solution of BIVP, 
except the right semi-plane is used. The only thing left is to use the linear change to return to the equation (1). \ $\blacksquare$

The same theorem can be formulated and proven in the left semi-plane.

\smallskip 
Let us describe the notion of uniqueness (global) of a solution of an initial value problem.

\smallskip 
{\bf Df.}\, 
Solution $y=\varphi(x)$ of an initial value problem, set at the point $(x_0,y_0)\in \widetilde G,$
is {\it unique on $\langle a,b\,\rangle \ni x_0,$} if $(x,\varphi(x))$ is a uniqueness point for any $x\in \langle a,b\,\rangle .$ 

\smallskip
The question arises: how to understand if the arc of the integral curve of the extension of the locally unique solution will not reach the non-uniqueness point? To answer this question, the way to constructively extract the sets, which consist of uniqueness points, is required.

\smallskip
{\bf Df.}\, Denote connected open set $G^o\subset G$ as  {\it a domain of uniqueness} of the equation~$(1),$
if every point of $G^o$ is a uniqueness point.
Denote set $\widetilde G^o=G^o\cup \widehat G^o,$ in which
$\widehat G^o$is a set of boundary points $G^o,$ which are also uniqueness points, as  {\it a uniqueness set.} 

\smallskip
It is implied by definition,  that the solution of the initial value problem is unique, until its graph is in the uniqueness set.

This means, that the definition of a domain of uniqueness $G^o$ coincides with the definition 1.2.1 from [8], 
where {\sl a connected open set $G^o\subset G$ is called a domain of uniqueness for equation (1), if for any point $(x_0,y_0)\in G^o,$ 
it is possible to choose such neighbourhood $V_\delta(x_0,y_0),$ that an arcs of any two integral curves, 
which pass through the point $(x_0,y_0)$ and are located in the neighbourhood $V_\delta(x_0,y_0),$ coincide. } 

Moreover, the latter definition can be automatically generalized on uniqueness set $\widetilde G^o.$

\smallskip
It is obvious, that the right-hand side of the equation (1) must satisfy additional require- ments, other than continuity, on $G^o\,(\widetilde G^o).$ The weaker requirements are, while still being sufficient, the stronger uniqueness theorem can be obtained. \goodbreak

\smallskip
{\bf Theorem} (on uniqueness set).\, 
{\it Assume function $f(x,y)$ in equation $(1)$ is defined and continuous on $\widetilde G.$ This function also satisfy Lipschitz condition with respect to $y$ locally on $\widetilde G^o=G^o\cup \widehat G^o,$ where the connected open set $G^o\subset G,$ and the set $\widehat G^o\subset \partial G^o\cap \widehat G.$ Then $\widetilde G^o$ is a uniqueness set of the equation $(1).$ }

P r o o f\,.\, Let us arbitrarily choose a point $(x_0,y_0)$ of $\widetilde G^o$ and prove that it is a uniqueness point.
Because $f\in {\rm Lip}_y^{loc}(\widetilde G^o),$ there are the closed $c$"=neighbourhood $\overline V_c(x_0,y_0)$ and the constant $L>0,$ 
such that  $f\in {\rm Lip}_y^{gl}(U_c)$ with the constant $L,$ where $U_c=\widetilde G^o\cap \overline V_c(x_0,y_0).$

If $(x_0,y_0)\in G^o,$ then such $c>0$ exists, that $U_c=\overline V_c(x_0,y_0),$ the solution of IIVP\!$(x_0,y_0)$ exists on some interval $(a,b)$ and, for any solution of this initial value problem, it can be achieved, restricting $(a,b)$ if needed, that its graph is in $U_c.$ 

Assume $(x_0,y_0)\in \widehat G^o.$ If there is no solution of BIVP\!$(x_0,y_0),$ then $(x_0,y_0)$ is a uniqueness point of the equation (1) by definition. If solution exists on some interval $\langle a,b\,\rangle ,$ such that $x_0\in \langle a,b\rangle \subset [x_0-c,x_0+c],$ 
then it can be achieved, restricting the interval if needed, that its graph is in $U_c.$ 

{\small In fact, it is obvious that, by restricting  $\langle a,b\,\rangle ,$ 
it can be achieved that the graph of the solution is in $c$"=neighbourhood $\overline V_c(x_0,y_0).$ 
Situation, when for $x<x_0$ and (or) $x>x_0$ graph, while being in $\widetilde G,$ is not in$\widetilde G^o,$ can be overcome by choosing such constant $c_1>c,$ that the global Lipschitz condition with the constant $L_1=L+1$ holds in the closed neighbourhood  $\overline V_{c_1}(x_0,y_0).$ 
As a result, taking the continuity of the function $f(x,y)$ into account, a connected open set $\widetilde G^o$ will expand, covering the arc of the integral curve in the small neighbourhood of the point $(x_0,y_0).$  }

\smallskip
The further proof is standard. It includes the writing down of two any solution in the integral form, using Lipschitz condition and Gronwall lemma.
\ $\blacksquare$ 

\smallskip
The special case of this theorem for a connected open sets is widely known.

\smallskip
{\bf Theorem} (on uniqueness in a connected open set).\, {\it 
Assume that in the equation $(1),$ the function $f(x,y)\in C(G), f\in {\rm Lip}_y^{loc}(G^o),$ the connected open set $G^o\subset G,$ 
then $G^o$ is a domain of uniqueness. }

\smallskip
In example\;1, a partial derivative of the right-hand side is not defined on the abscissa axis, and the integral curve 
of the singular boundary solution $y(x)\equiv 0,\ x\in [0,+\infty)$ is not in the uniqueness set $\widetilde G^o=\{(x,y)\!:\,x\ge 0,\ y>0\}.$

\smallskip
In practice, checking that the function $f(x,y)$ satisfies the local Lipschitz condition is harder, than verifying that its partial derivative 
with respect to $y$ exists. Therefore, it is convenient to use the following proposition to study uniqueness of a solution. 

\smallskip
{\bf Theorem} (on uniqueness set; weak).\, {\it 
Assume, that in the equation $(1)$ the function $f(x,y)$ is continuous on $\widetilde G,$ the function $\partial f(x,y)/\partial y$ 
is continuous on $G^o\subset G.$ Then the set $\widetilde G^o=G^o\cup \widehat G^o,$ 
where $\widehat G^o\subset \partial G^o\cap \widehat G$ and consists of points, in which $\partial f(x,y)/\partial y$
can be continuously redefined, is a uniqueness set, if for any point $(x_0,y_0)\in \widehat G^o,$ 
such $c$-neighbourhood $\overline V_c(x_0,y_0)$ exists, that the set $\widetilde G^o\cap \overline V_c(x_0,y_0)$ is convex with respect to $y.$ }

P r o o f\,. Consider arbitrarily chosen point $(x_0,y_0)\in \widetilde G^o.$

Because the function $\partial f(x,y)/\partial y$ is continuous at $(x_0,y_0),$ such $\delta$ exists, that $0<\delta\le c,$ 
where $c$ is given in theorem's formulation, and for any point $(x,y)\in U_\delta=\widetilde G^o\cap \overline V_\delta(x_0,y_0),$ 
the inequation $|\partial f(x,y)/\partial y-\partial f(x_0,y_0)/\partial y|\le 1$ is valid. 
Therefore $U_\delta$ is convex with respect to $y$
and $|\partial f(x,y)/\partial y|\le L=|\partial f(x_0,y_0)/\partial y|+1$ for any point $(x,y)\in U_\delta.$

According to Lagrange theorem, for any two points $(x,\hat y),(x,\tilde y)\in U_\delta,$
$\hat y<\tilde y,$ such $y^*(x)\in(\hat y,\tilde y)$ exists, that
$\displaystyle f(x,\tilde y)-f(x,\hat y)=f_y'(x,y^*(x))(\tilde y-\hat y).$

Here, the point $(x,y^*(x))\in U_\delta,$ because the set $U_\delta$ is convex with respect to $y.$
Therefore, the inequation $|f(x,\tilde y)-f(x,\hat y)|\le L|\tilde y-\hat y|$ is valid in $U_\delta,$ 
which means that $f\in {\rm Lip}_y^{gl}(U_\delta).$
Then$f\in {\rm Lip}_y^{loc}(\widetilde G^o),$ and according to the theorem on uniqueness, $\widetilde G^o$ is a uniqueness set. $\blacksquare$

\smallskip
{\bf Example 3.}\, Consider equation 
$$y'=\big\{2|y|^{3/2}\ \,\hbox{when}\ \,x\le 0,\ 2(|y|-x^2)^{3/2}+x\cdot {\rm sign}\,y\ \;\hbox{when}\ \,x\ge 0\big\},\eqno (8)$$ 
its right-hand side $f(x,y)$ and its partial derivative with respect to $y$ are continuous on connected set 
$\widetilde G=\big\{(x,y)\!:\,y\in \mathbb{R}^1\ \hbox{when}\ x\le 0,$ $|y|-x^2\ge 0\ \hbox{when}\ x\ge 0\big\},$ 
and $\widehat G=\{y=\pm x^2\ \hbox{when}\ x\ge 0\},$ 
functions $y=\pm x^2,$ $x\in [0,+\infty)$ are complete boundary solutions. 

The uniqueness set $\widetilde G^o$ in equation (8) is $\widetilde G\backslash \{(0,0)\},$
because the origin is a non-uniqueness point due to the integral curves of boundary solutions having a common tangent at that point.

Other than the interior particular solution $y(x)=0,\ x<0$ and the mixed boundary solutions $y=\pm x^2,\ x\ge 0,$ 
a general solutions of the equation (8) can be written down in upper and lower semi-planes correspondingly:

\,$y=\bigg [\begin{matrix} \big\{(C-x)^{-2},\, x\le 0;\ (C-x)^{-2}+x^2,\, x\in[0,C)\big\}\ \hbox{ where }\ C>0, \\ 
         (C-x)^{-2},\, x<C\ \hbox{ where }\ C\le 0 \hfill \end{matrix}$ 
				
\,$y=\bigg [\begin{matrix} -(x-C)^{-2}-x^2,\  x\in(C,+\infty)\ \hbox{ where }\ C\ge 0, \hfill \\ 
          \big\{-(x-C)^{-2},\, x\in(C,0];\ -(x-C)^{-2}-x^2,\, x\ge 0\big\}\ \hbox{ where }\ C<0. \end{matrix}\quad \square$  

\includegraphics[scale=0.53]{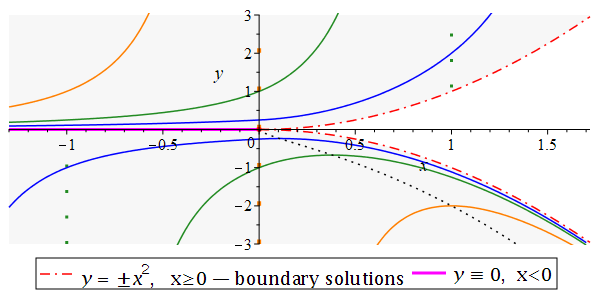} 

\smallskip 
Example \;3 demonstrates the necessity of the condition on local convexity with respect to $y$ of the uniqueness set in weak theorem on uniqueness set.
Moreover, as expected, for any neighbourhood of the point $(0,0)$ function $f$ does not satisfy the global Lipschitz condition, 
because, for any $x>0,$ the Lipschitz equation calculated at the points $(x,-x^2),(x,x^2)\in \widehat G$ takes form $x\le L x^2.$ Therefore, 
for any $L>0,$ such constant $c>0$ exists, that this inequality is not valid on the set $\widetilde G^o\cap \overline V_c(0,0),$ because $x\le c.$ 

Let us note in conclusion, that the condition of the existence of continuous partial derivative of the function $f(x,y)$ with respect to $y$ 
in formulation of the weak theorem on uniqueness set is not required. This fact is supplemented by counterexamples.

The partial derivatives with respect to $y$ of the right-hands sides of the equation (5) and (6) are not defined on abscissa axis, however, the graph of the origin interior solution of the equation (5) contain the countable number of uniqueness points, and all points of the graph of the origin boundary solutions of the equation (6), except the origin itself, are uniqueness points.
\newpage

\begin{center} REFERENCES \end{center}

{\small 
1. Basov V.V., Iljin Yu.A. УOn the Existence of a Solution to the Cauchy Initial Boundary Value ProblemФ, 
\textit{Vestnik St.\,Petersburg University, Mathematics}. Vol.\,53. No 2, (2020).

2. Basov V.V., Iljin Yu.A. УOn the Cauchy Problem Set on the Boundary of the Ordinary Differential EquationТs Domain of DefinitionФ, 
\textit{Vestnik St.\,Petersburg University, Mathematics}. Vol.\,53. No 4, (2020).

3. Lakshmikantham V., Leela S. Differential and integral inequalities; theory and applications, vol.I, Academic Press, New York, 1969.

4. Krasnosel'skii M.A. The operator of translation along the trajectories of differential equations. 
Translation of Mathematical Monographs., Vol. 19, Providence, R.I., 1968.

5. Bibikov Yu.N., Pliss V.A.,
УOn the Dependence of Initial Values of the Maximal Interval of Existence of a Solution of a Differential Equation Ф,
{\it Vestnik St. Petersburg University. Mathematics}, {\bf 47}(4), 141--144 (2014).

6. Basov V.V., {\it Ordinary differential equations: lectures and practical seminars}. 
(S.-Petersburg University Press, Saint-Petersburg, 2023) [In Russian].

7. Hartman Ph., {\it Ordinary differential equations} (John Willey and Sons, New York London Sydney, 1964).

8. Bibikov Yu.N., {\it General course of ordinary differential equations}. (S.-Petersburg University Press, Saint-Petersburg, 2005) [In Russian].

}

\medskip
{\footnotesize
Author's information:

{\it Vladimir V. Basov} --- vlvlbasov@rambler.ru }

\end{document}